\documentclass[amsfonts,12pt]{article}
 \usepackage{latexsym,epsfig,graphpap,graphics,amsmath}
\usepackage{amssymb,amsmath}
\newtheorem{theorem}{Theorem}
\def \curl{\overrightarrow{\mb{curl}}}

\newtheorem{corollary}[theorem]{Corollary}

\newtheorem{definition}[theorem]{Definition}

\newtheorem{lemma}[theorem]{Lemma}

\newtheorem{proposition}[theorem]{Proposition}

\def \var{\varphi}

\def \h{\hat}
\def \p{\partial}
\def \u{\vec{u}}
\def \f{\vec{f}}

\def \non{\nonumber\\[2mm]}

\def \n1{\newpage}

\def \no1{\noindent}

\def \fr{\frac}
\def \ti{\times}

\def \seq{\subseteq}
\newcounter{proposition}
\newcounter{definition}

\def \tu{\vartriangle}
\def \td{\triangledown}

\def \h{\hat}
\def \bl{\begin{lemma}}
\def \bcor{\begin{corollary}}
\def \ecor{\end{corollary}}
\def \el{\end{lemma}}
\def \1{\begin{eqnarray}}
\def \2{\end{eqnarray}}
\def \3{\begin{eqnarray*}}
\def \4{\end{eqnarray*}}
\def \6{\vspace*{7mm}}
\def \bqn{\begin{equation}}
\def \eqn{\end{equation}}

\def \s1{\subseteq}

\def \mb{\mbox}

\def \bt{\begin{tabular}}
\def \et{\end{tabular}}

\def \l{\left}

\def \r{\right}
\def \hs1{\hspace*{3mm}}
\def \q2{\hspace*{9mm}}

\def \un1{\underline}
\def \mb{\mbox}
\def \vs1{\vspace{4mm}}

\def \ba{\begin{array}}
\def \ea{\end{array}}
\def \Om{\Omega}
\newcommand{\ec}{\end{center}}
\newcommand{\bc}{\begin{center}}
\newcommand{\be}{\begin{equation}}
\newcommand{\ds}{\displaystyle}

\newcommand{\ee}{\end{equation}}
\newcommand{\bn}{\begin{enumerate}}
\newcommand{\en}{\end{enumerate}}
\newcommand{\bi}{\begin{itemize}}
\newcommand{\ei}{\end{itemize}}
\begin{document}
\title{Finite Element Technique for Solving the Stream
Function Form of a Linearized Navier-Stokes Equations
Using Argyris Element}
\author{\small  F. Fairag and N. Almulla}
\date{}
\maketitle

\begin{abstract}
The numerical implementation of finite element discretization method
for the stream function formulation of a linearized Navier-Stokes
equations is considered.  Algorithm 1 is  applied
                        using Argyris element.  Three global
orderings of nodes are selected and registered in order to conclude
the best banded structure of matrix and a fluid flow calculation
is considered to test a problem which has a known solution.
Visualization of global node orderings,
matrix sparsity patterns and stream function contours
                 are displayed showing the main features of the flow.
\end{abstract}

\vs1


\vs1

\section{Introduction}
The numerical treatment of nonlinear problems that
arise in areas such as fluid mechanics often requires solving
large systems of nonlinear equations.  Many methods have been
proposed that attempt to solve these systems efficiently;
one such class of methods is the finite element method, the most widely
used technique for engineering design and analysis.

The Navier-Stokes equations may be solved using either the primitive
variable or stream function formulation.  Here, we use the stream
function formulation.  The attractions of the stream function
formulation are that the incompressibility constraint is automatically
satisfied, the pressure is not present in the weak form, and there is
only one scalar unknown to be determined.
            The standard weak formulation
of the stream function first appeared in 1979 in [9];
a general
analysis of convergence for this formulation has been done
                 in [4, 5].

The goal of this paper is to demonstrate that the method can be
implemented to approximate solutions for incompressible viscous
flow problems.
\section{Governing Equations}
\renewcommand{\theequation}{\thesection.\arabic{equation}}
\setcounter{equation}{0}
Consider the Navier-Stokes equations describing the flow
of an incompressible fluid
\1
- \mb{Re}^{-1} \tu \u + (\u \cdot \td)\u +
\td p & = &  \f ~~ \mb{in} ~~ \Om, \\
\td \cdot \u & = &    0 ~~ \mb{in} ~~ \Om, \\
           \u & = &    0 ~~ \mb{on} ~~ \p \Om, \2
where $\u = (u_1, u_2)$ and $p$ denotes the unknown velocity
and pressure field, respectively, in a bounded, simply connected
polygonal domain $ \Om \seq R^2$; $\f$ is a given body force; and
$\mb{Re}$ is the Renolds number.

The introduction of a stream function $\psi(x,y)$ defined
by
\[ u_1 = - \fr{\p \psi}{\p y}, ~ u_2 =
\fr{\p \psi}{\p x}\]
means that the continuity equation (2.2) is satisfied identically.
The pressure may then be eliminated from (2.1) to give
\1
\mb{Re}^{-1} \tu^2 \psi - \psi_y \tu \psi_x +
\psi_x \tu \psi_y & = &
\overrightarrow{\mb{curl}}
\f ~~\mb{ in }
\Om \\
\psi & = & 0 ~~ ~~\mb{on} ~~ \p \Om \\
\fr{\p \psi}{\p \h{n}} & = & 0            ~~~~ \mb{on} ~~ \p \Om, \2
where $\h{n}$ represents the outward unit normal to $\Om$.
Equation (2.4) is a nonlinear fourth-order partial
differential equation which turns into the known linear biharmonic
fourth-order partial differential equation by omitting the second and
third terms in the left-hand side of the equation.
In order to write (2.4)\,--\,(2.6) in a
variational form, we define the Sobolev spaces
\1
H^1 (\Om) & = & \{v:v \in L^2(\Om), Dv\in L^2(\Om)\},\\
H^2_0 (\Om) & = & \{v:v \in H^1(\Om),  v = 0 ~~ \mb{on} ~~ \p \Om \},\\
H^2   (\Om) & = & \{v:v \in L^2(\Om), Dv\in L^2(\Om), D^2v\in L^2\}, \\
H^2_0 (\Om) & = & \{v:v \in H^2(\Om): v = \fr{\p v}{\p n} = 0, ~~
\mb{on} ~~ \p \Om\},
\2
where $L^2(\Om)$ is the space of square integrable functions
on $\Om$ and $D$ represents differentiation with respect to $x$ or $y$.
For each $\var\in H^1(\Om)$, define $\overrightarrow{\mb{curl}}
~ \var = \l(\ba{c}\var_y \\ -\var_x\ea \r)$. The
following linearized weak form of equations (2.4)\,--\,(2.6) is
considered:
\1
&& \mb{Find } \psi\in H^2_0(\Om) \mb{ such that for all }
\var\in H^2_0 (\Om), \non
&&
a(\psi, \var) + b(\psi^*; \psi, \var) = \ell(\var),
\2
where
\3
   a(\psi, \var) & = &
\mb{Re}^{-1} \int_\Om \tu \psi \;\cdot \tu \var d\Om, \\
b(\xi;\psi, \var) & = & \int_\Om \tu \xi (\psi_y
\var_x - \psi_x \var_y)d\Om, \\
\ell(
\var) & = & (\f, \curl ~ \var) = \int_\Om \f \cdot
\curl ~ \var d\Om,
\4
where $\psi^*$ is a fixed given function (a primitive
approximation for $\psi$).
\section{Finite Element Discretization}
\setcounter{equation}{0}
For the standard finite element discretization of (2.11), we choose
conforming finite element subspace $X \subset H^2_0(\Om)$.
We then seek $\psi^h \in X$ such that for all $\var^h\in X$,
\bqn
a(\psi^h, \var^h) + b(\psi^{* h}, \psi^h, \var^h) = \ell(\var^h). \end{equation}
          The existence and uniqueness for the
solution of the discrete problem (3.1) has been proved in [1].
          Once the finite element
spaces are prescribed, problem (3.1)
reduces
to solving a system of algebraic equations.  Various
iterative methods can be used to solve
           problem (3.1).
           In this paper the following algorithm
           has been applied to solve problem (3.1) for
           a fixed mesh of size $h$.
In each iteration in Algorithm 1, we need to solve a linear system. The
    resulting linear system is nonsymmetric whose symmetric
part is positive definite. Moreover, the resulting matrix is a sparse
matrix.
We choose the  conjugate gradient stabilized (BICG\,STAB) method,
which requires two matrix-vector products and four inner products in
each
iteration.
                                                              \\[4mm]
 \noindent\begin{tabular}{|l|} \hline
\noindent {\bf Algorithm 1} (the finite element algorithm)
\\[2mm]
       Given
Max-iteration \& Tolerance\\
       Given $\psi_0$ as a starting guess\\
       For $i=1:$ Max-iteration \\
       \qquad Solve the linear system on the mesh $X\subset H^2_0$ for
$\psi_i$: \\
\hspace*{2cm} $\ds a(\psi_i, \var) + b(\psi_{i-1}, \psi_i, \var) = \ell(\var)
~~ \forall\;
\var \in X$ \\
       \qquad If $\|\psi_i - \psi_{i-1}\| \leq$ Tolerance \& Residual $\leq $
Tolerance \\
       \qquad Stop\\
End \\ \hline \et
\section{Finite Element Space}
\setcounter{equation}{0}
The inclusion $X \subset H^2_0(\Om)$ requires the use of finite
element functions that are continuously differentiable over $\Om$.
We choose Argyris triangle as a finite element space for the stream
function formulation. We will impose boundary conditions by setting all
the degrees of freedom at the boundary nodes to be zero and the normal
derivative equal to zero at all vertices and nodes on the boundary.
In Argyris triangle, the functions are quintic polynomials within each
triangle and the 21 degrees of freedom are chosen to be the function
value, the first and second derivatives at the vertices, and the normal
derivative at the midsides.
\section{Global Node Orderings}
\setcounter{equation}{0}
For linear systems derived from a partial differential equation, each
unknown corresponds to a node in the discretization mesh.  Different
orderings of unknowns correspond to permutations of the coefficient
matrix.  Since the convergence speed of iterative method may depend on the
ordering used,  the next three options have been
considered.
\begin{enumerate}
\item Start numbering the nodes at the
vertices with six successive numbers at each node.  Then
the numbering of the midpoint nodes is performed starting with
the horizontal sides, followed by the vertical and oblique
side.
This is illustrated in Figure (5.1).
\item  Start numbering the
function values of the nodes at vertices followed by the
derivatives, and then go back to the midpoint nodes
using the same ordering as above.
           This is illustrated in Figure (5.2).
\item Start numbering
the nodes at the vertices with six successive numbers by
skipping a node alternatively, and then go back to the
midpoint nodes of horizontal, vertical and oblique sides
using         the same orderings as above.
                This is illustrated in Figure (5.3).
\end{enumerate}
\begin{figure}[h] \centering
\epsfig{file=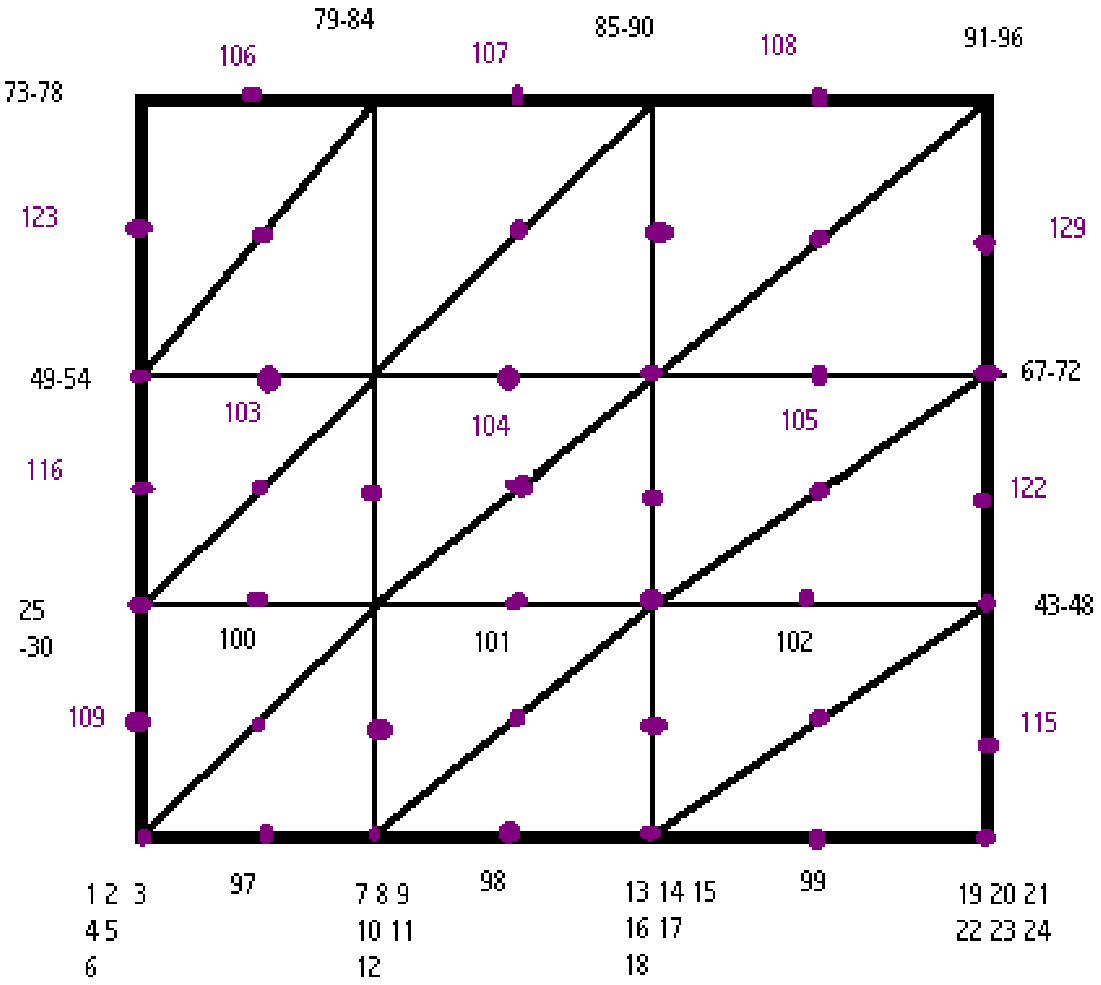,
width=8.5cm,height=8.5cm} \\[-2mm]
Figure 5.1
\end{figure}
\begin{figure}[p]\centering
\epsfig{file=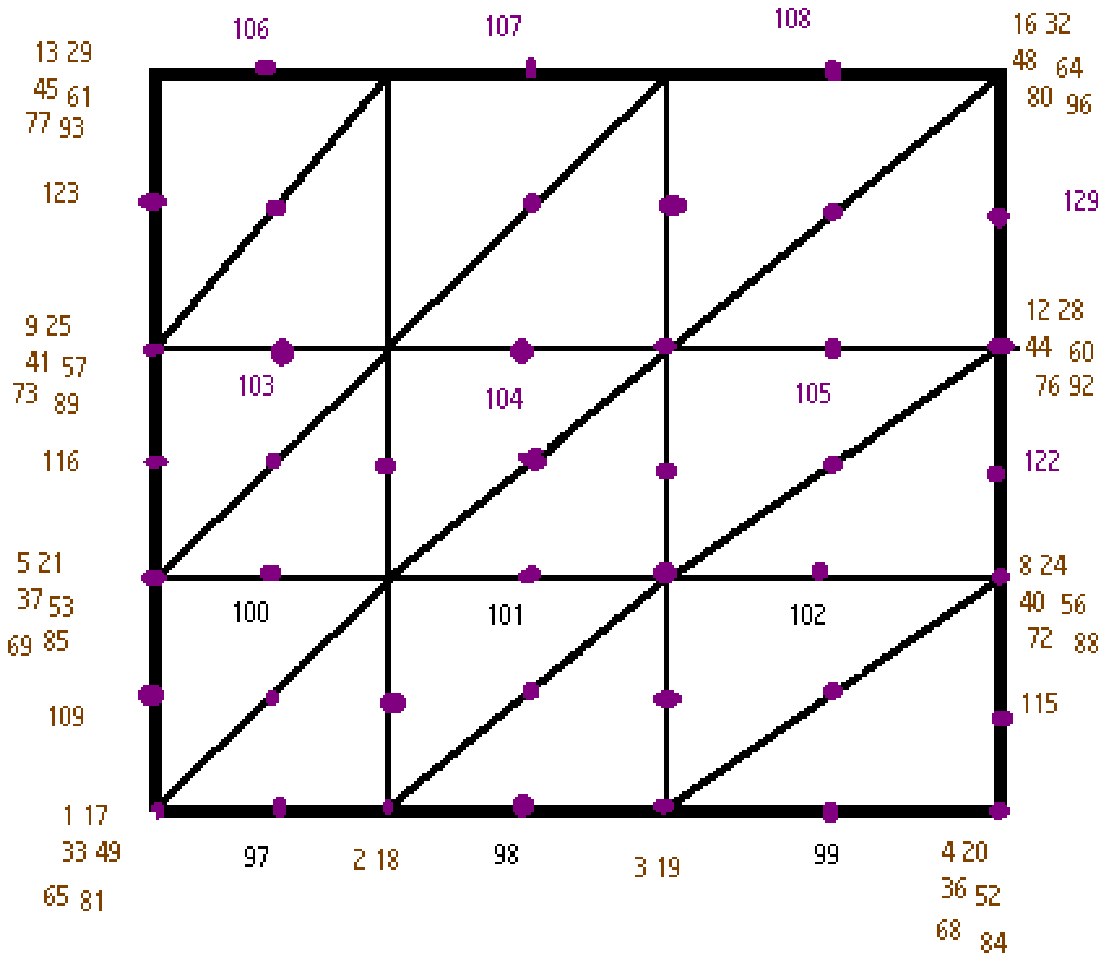,
width=9.75cm,height=9.75cm} \\[-2mm]
Figure 5.2
\end{figure}
\begin{figure}[p]\centering
\epsfig{file=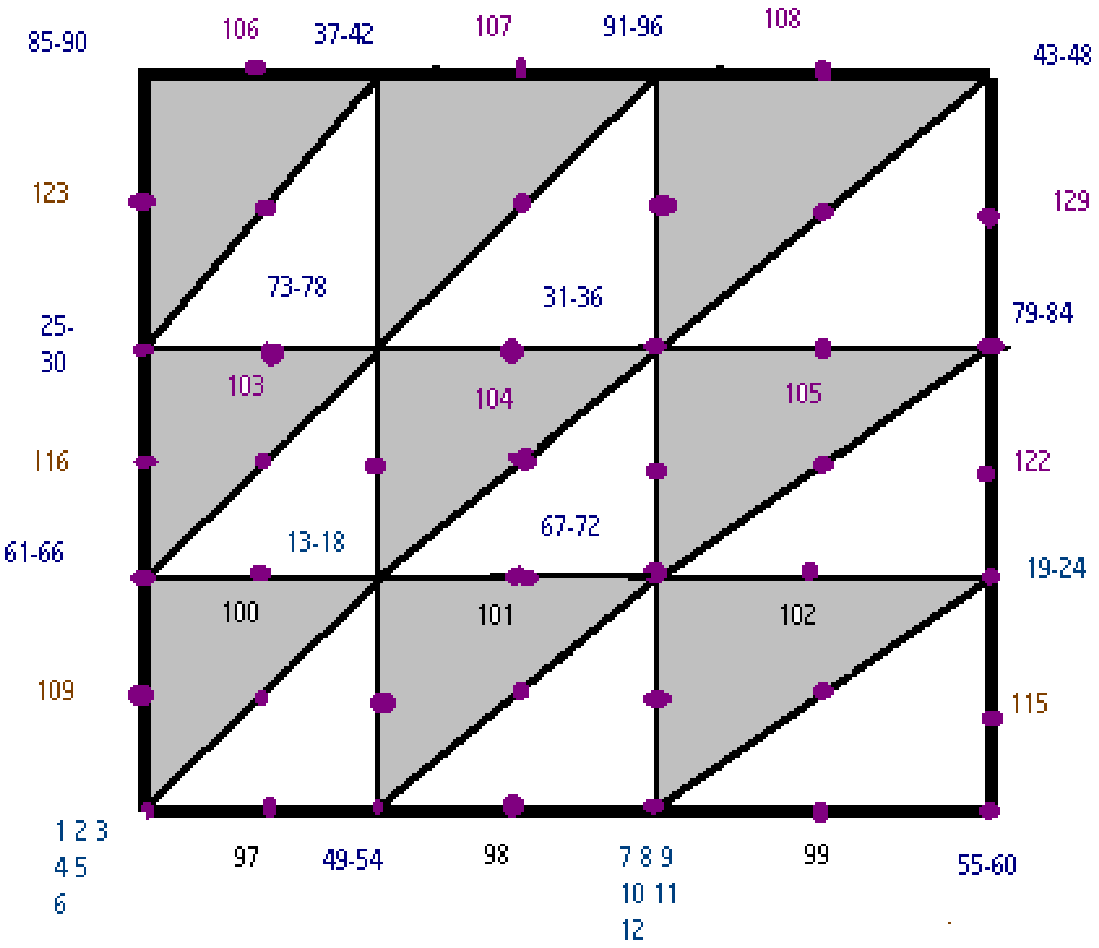,
width=9.75cm,height=9.75cm} \\[-2mm]
Figure 5.3
\end{figure}
\clearpage
It can be concluded that the first option is the most appropriate
since it has the best matrix property, being a banded structure, as
shown in Figures 5.4, 5.5, 5.6
 which visualize, respectively, the location of the
 nonzero elements for the three selected global orderings.
Therefore this ordering has been used in the next section.

\vspace{6mm}

{\begin{figure}[h]\centering
\epsfig{file=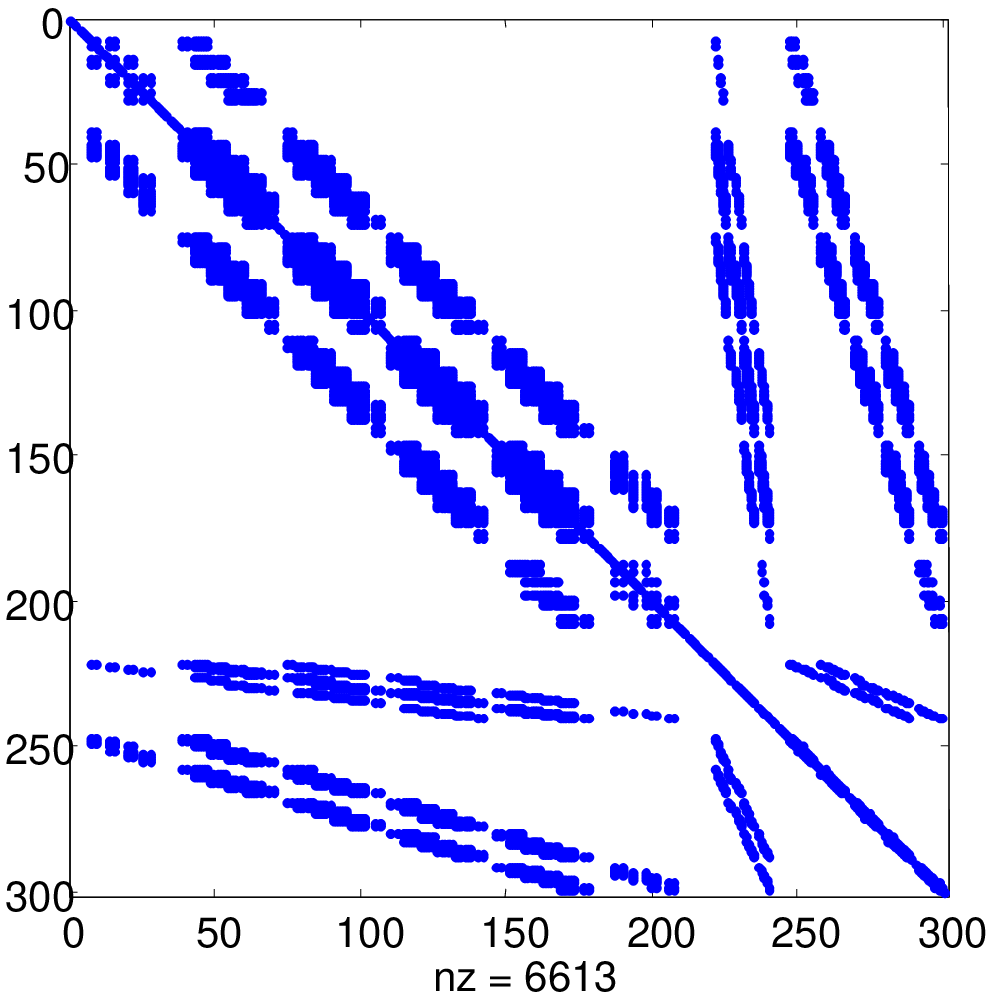, width=7.75cm,height=7.75cm} \\[2mm]
Figure 5.4
\end{figure}
\begin{figure}[h,b]\centering
\epsfig{file=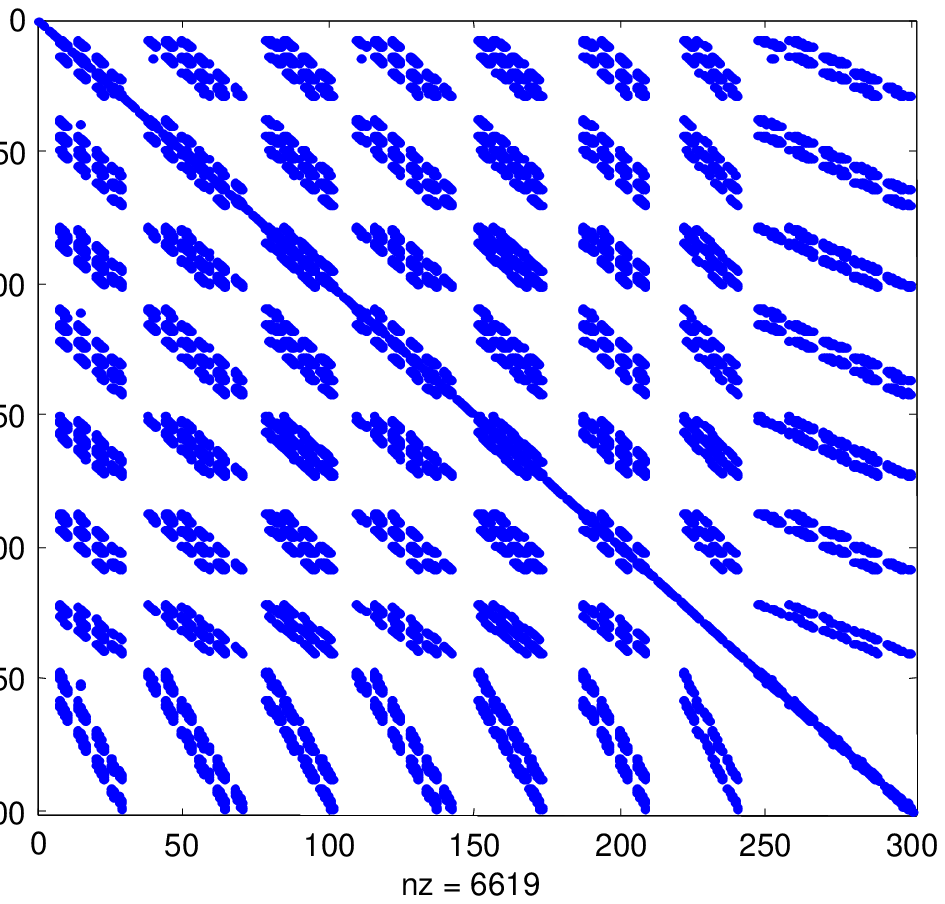,
 width=7.75cm,height=7.75cm
} \\[2mm]
Figure 5.5
\end{figure}
\clearpage
\begin{figure}[t]\centering
\epsfig{file=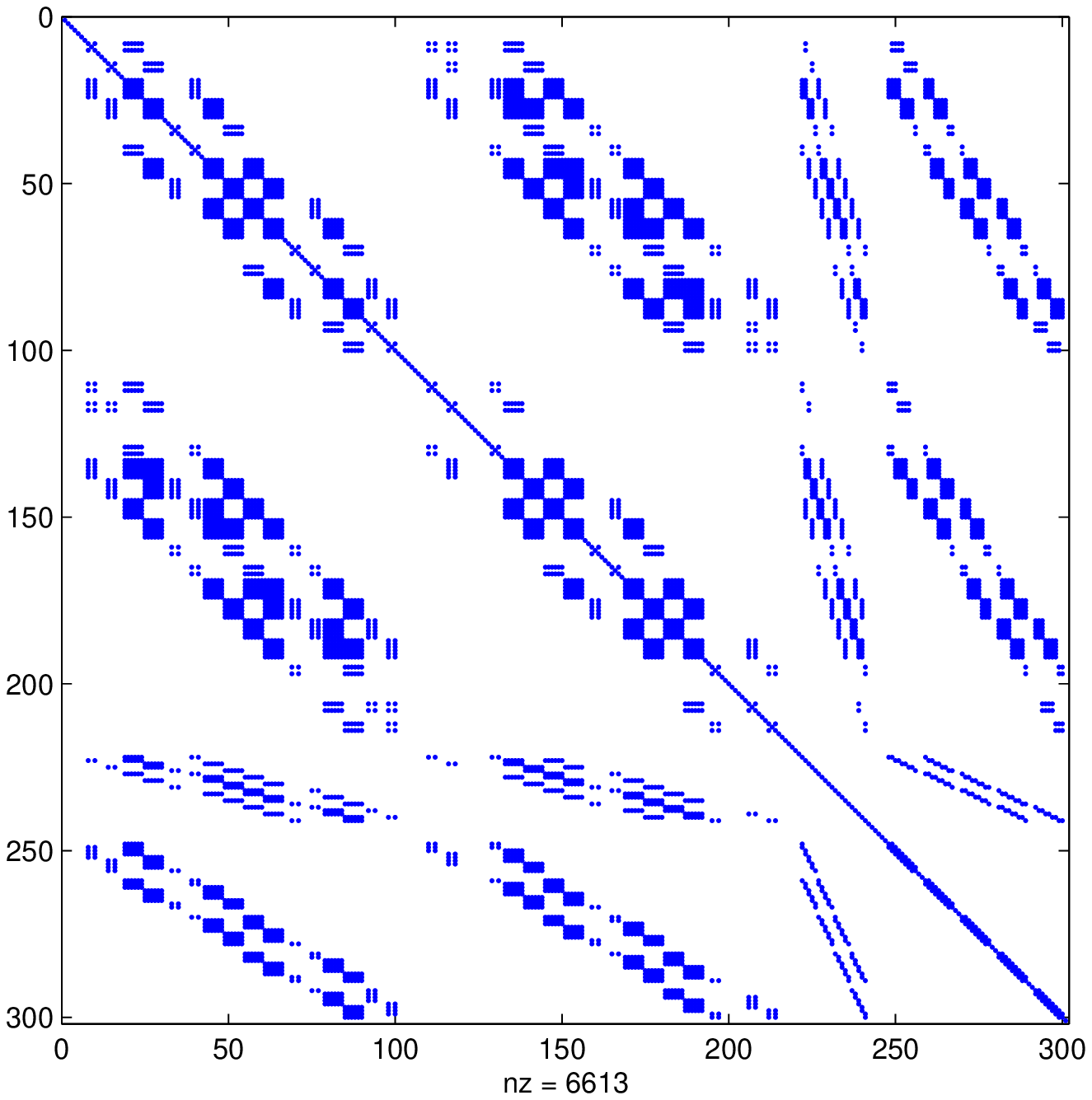,
 width=9.25cm,height=9.25cm
} \\[2mm]
Figure 5.6
\end{figure} }

\vspace{4mm}

Furthermore, the following table illustrates the cpu-time,
the number of arithmetic operations, and the number of bicgstab
iterations which have been computed by performing the
           suggested global numberings in a specified problem:
\bc
\bt{|l|c|c|c|} \hline
& cpu-time & n.c.o. & bicgstab iter. \\ \hline
1. & 28.45 & 227800974 & 282.5 \\ \hline
2. & 29.77 & 231648608 & 289.0\\ \hline
3. & 29.06 & 239406016 & 297.5 \\ \hline
\et
\ec
\section{Numerical Example}
In this section we describe some numerical results obtained
by implementing the finite element algorithm, for which we have
an exact solution. The region $\Om$ is the
unit square $\{0 < x < 1, ~ 0 < y < 1\}$ and for the finite
element discretization, we use the Argyris elements.
The biharmonic equation along with equation (2.11) was
solved in the same way.
                                                     They were solved
until
the norm of the difference in successive iterates and the
norm of the residual were within a fixed tolerance.

We consider  as a test example the two-dimensional Navier-Stokes
equations (2.1)\,--\,(2.3) on the unit square $\Om = (0, 1)^2$, where
we define the right-hand side by $f: = - \mb{Re}^{-1}
\tu \u + (\u \cdot \td)\u + \td p$ with
the following prescribed exact solution:
\3
u   =   \l(\ba{c}\psi_y \\ -\psi_x\ea \r) \mb{ with } ~~
\psi(x,y) & = & x^2 (x-1)^2 y^2 (y-1)^2, \\
p & = & x^3 + y^3 - 0.5.
\4
For this test problem, all requirements of the theory concerning
the geometry of the domain and the smoothness of the data are
satisfied. Moreover, the stream function $\psi(x, y)$
satisfies the
boundary conditions of the stream function equation
of the Navier-Stokes equations.

In all numerical calculations in this example, we have
used the Argyris elements with $\mb{Re } = 1$ and $\mb{tol } =
10^{-5}$.
We pick three values of $h$: 1/3, 1/5, and 1/9. The cpu-time,
the number of PCG iterations, and
the error $|\psi_c - \psi_e|$
for different values of $h$
are tabulated in Table 6.1 based
on four quadrature points, and Table 6.2 based on six quadrature
points
            for the biharmonic equation.
           The  cpu-time, the number of BICGSTAB iterations, errors, and
the number of the mathematical operations are
tabulated in Table 6.3 based on six quadrature points for
    the linearized equation (2.11), where
$\psi^*$ is computed by solving the biharmonic equation
as initial guess, and where
\3
\renewcommand{\arraystretch}{1.25}
\mb{n.q.p.} & = & \mb{number of quadrature points,} \\
\mb{n.c.o.} & = & \mb{number of operations,} \\
\mb{pcg-itr.} & = & \mb{number of PCG iterations,} \\
\mb{bicgstab itr.} & = & \mb{number of BICGSTAB iterations,} \\
\psi_c & = & \mb{The computed solution of the function values,} \\
\psi_e & = & \mb{The exact solution of the function values.}
\4

\[
\ba{|c|c|c|c|c|c|} \hline
h & \mb{n.q.p} & \mb{cpu time} & \mb{n.c.o.} & \mb{Error} &
\mb{pcg-itr} \\
&&&& = \psi_c - \psi_e & \\ \hline
1/3  & 4 & 1.3200 & 17601719 & 0.1644\ti 10^{-3} & 72 \\ \hline
1/5  & 4 & 6.7500 & 246112068 & 0.1899\ti 10^{-3} & 211 \\ \hline
1/9  & 4 & 64.4300 & 3.7401\ti 10^9 & 0.1376\ti 10^{-3} & 437 \\ \hline
\ea \]
\bc Table 6.1 \ec
\[
\ba{|c|c|c|c|c|c|} \hline
h & \mb{n.q.p} & \mb{cpu time} & \mb{n.c.o.} & \mb{Error} &
\mb{pcg-itr} \\
&&&& = \psi_c - \psi_e & \\ \hline
1/3  & 6 & 1.8200 & 16476551 & 0.1984\ti 10^{-3} & 72 \\ \hline
1/5  & 6 & 7.4700 & 249166334 & 0.1935\ti 10^{-3} & 209 \\ \hline
1/9  & 6 & 64.7100 & 3.6482\ti 10^9 & 0.1379\ti 10^{-3} & 454 \\ \hline
\ea \]
\nopagebreak[4]
\bc Table 6.2  \ec
\[
\arraycolsep=.05in
\ba{|c|c|c|c|c|c|c|} \hline
h & \mb{cpu time} &  \mb{n.c.o.}
&  |\psi_c - \psi_e|_0
&  |\psi_c - \psi_e|_1
&  |\psi_c - \psi_e|_2 & \mb{bicgstab itr.} \\ \hline
1/3 & 13.24 & 17528466 & 0.2589 \ti 10^{-3} & 0.1294 \ti 10^{-2} &
0.1692 \ti 10^{-1} & 100.5 \\ \hline
1/5 & 35.32 & 227800974 & 0.2148 \ti 10^{-3} & 0.1062 \ti 10^{-2} &
0.1048 \ti 10^{-1} & 282.5 \\ \hline
1/9 &179.88 & 3.50366282\ti 10^9  & 0.1423\ti  10^{-3} & 0.6986 \ti 10^{-3} &
0.6016 \ti 10^{-2} & 567.5 \\ \hline
\ea \]
\bc Table 6.3 \ec
\begin{figure}[h,b]\centering
\epsfig{file=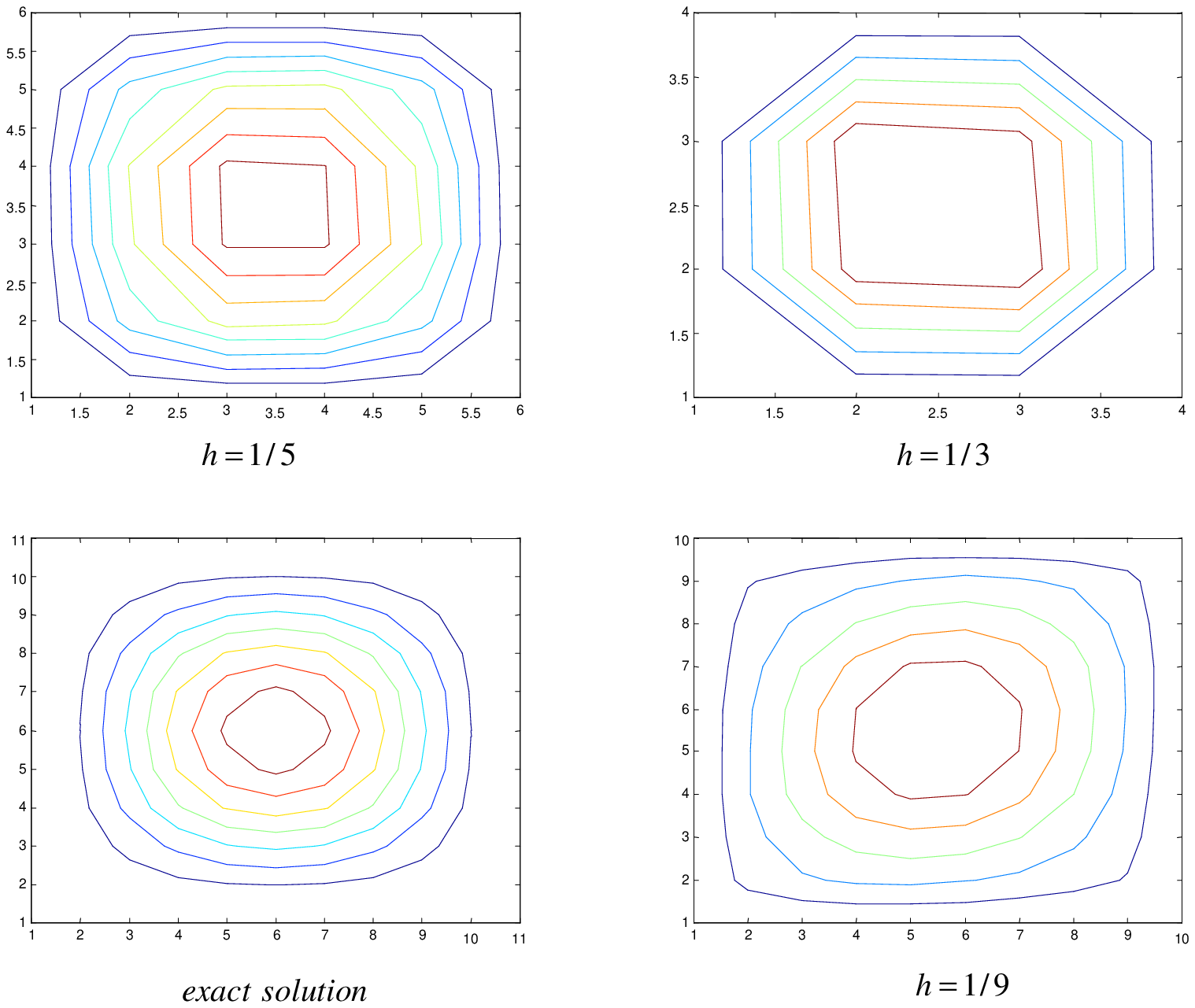}\\[2mm] 
Fig. 6.1: Streamlines for $\ds h = 1/3, 1/5, 1/9$ with $\mb{Re }= 1$
using finite element method on the linear problem with four
quadrature points as shown in Table 6.1
\end{figure}
\begin{figure}\centering
\epsfig{file=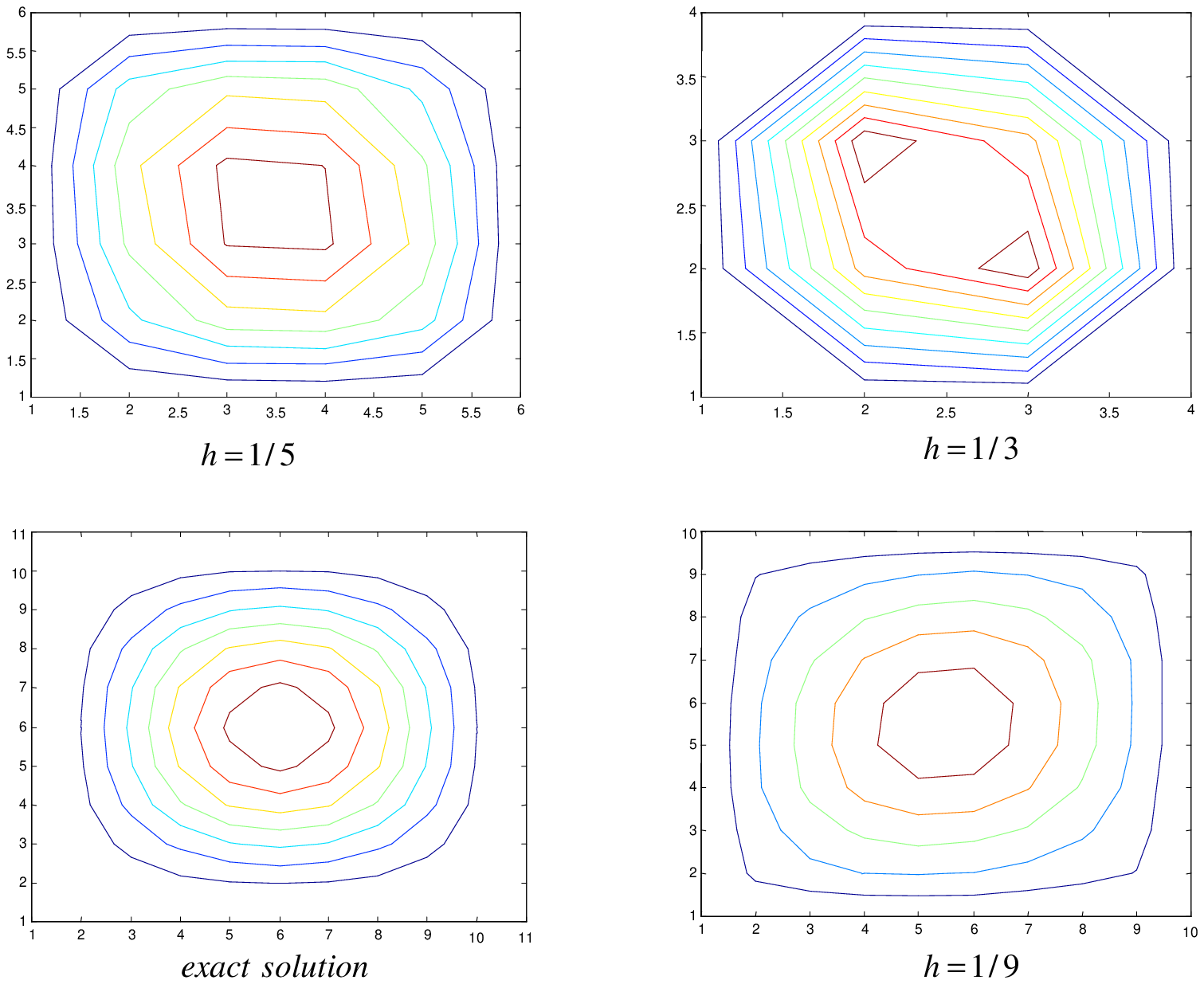}\\[2mm] 
Fig. 6.2: Streamlines for $\ds h = 1/3, 1/5, 1/9$ with $\mb{Re }= 1$
using finite element method on the linear problem with six
quadrature points as shown in Table 6.2
\end{figure}
\begin{figure}\centering
\epsfig{file=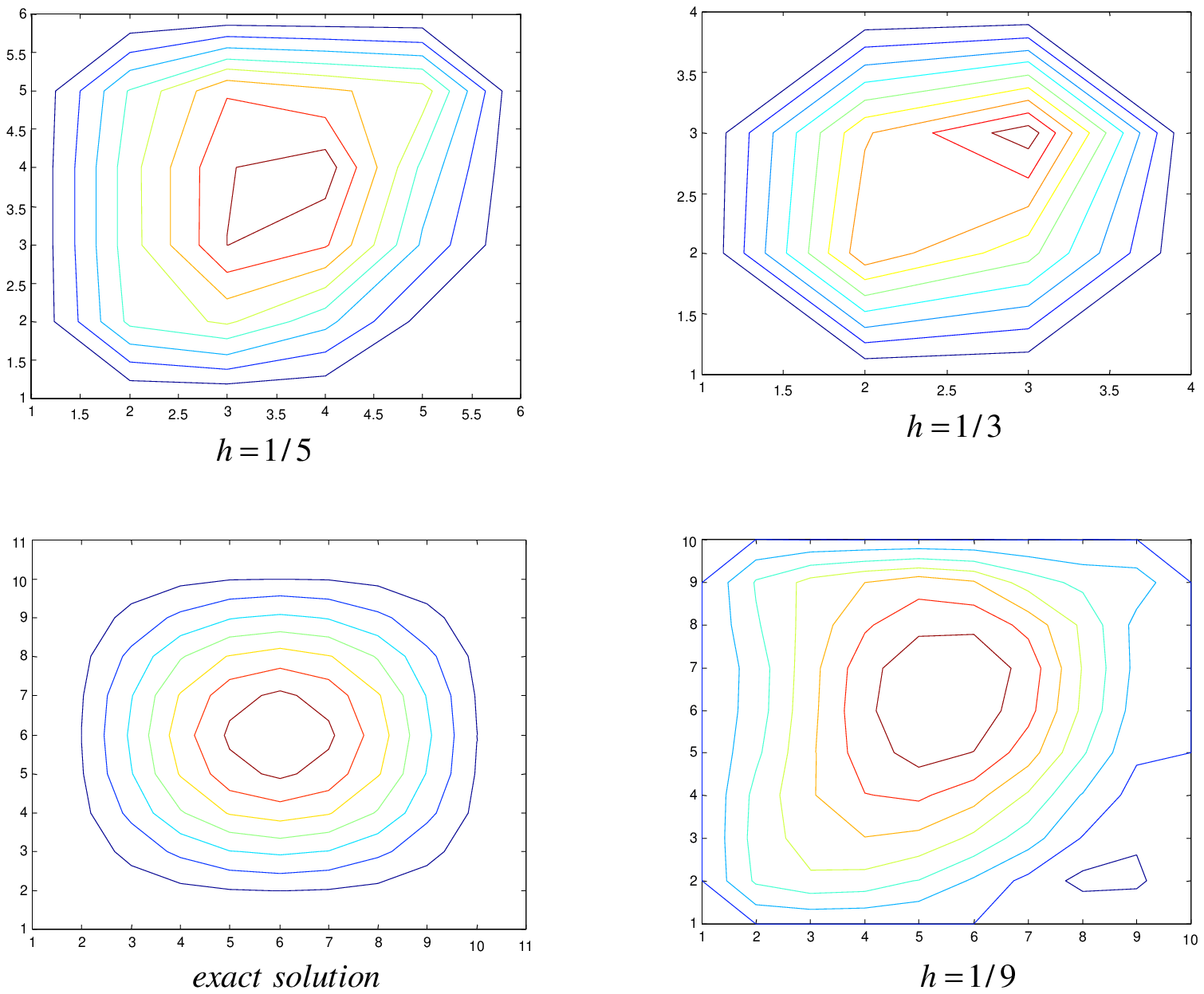}\\[2mm]
Fig. 6.3: Streamlines for $\ds h = 1/3, 1/5, 1/9$ with $\mb{Re }= 1$
using finite element method on the linear problem with six
quadrature points as shown in Table 6.3
\end{figure}

\end{document}